\documentclass[11pt]{article}
\usepackage{amsmath}
\usepackage{dsfont}
\usepackage{mathrsfs}
\usepackage{amsmath,amssymb}
\usepackage{amsfonts}
\usepackage{hyperref}
\usepackage{amsthm}
\usepackage{graphicx}
\usepackage{subfigure}
\usepackage{xcolor}
\renewcommand{\qed}{\hfill\small{$\square$}\normalsize}

\hfuzz=\maxdimen
\theoremstyle{definition}
\newtheorem{lemma}{Lemma}[section]
\newtheorem{definition}[lemma]{Definition}
\newtheorem{proposition}[lemma]{Proposition}
\newtheorem{theorem}[lemma]{Theorem}
\newtheorem{corollary}[lemma]{Corollary}

\newtheorem{remark}{Remark}

\newtheorem{claim}{Claim}

\numberwithin{equation}{section}

\renewcommand{\qed}{\hfill\small{$\square$}\normalsize}

\DeclareFixedFont{\Acknowledgment}{OT1}{cmr}{bx}{n}{14pt}
\begin{document}

\title{On the deformation of inversive distance circle packings, III}
\author{Huabin Ge, Wenshuai Jiang}
\maketitle

\begin{abstract}
Given a triangulated surface $M$, we use Ge-Xu's $\alpha$-flow \cite{Ge-Xu1} to deform any initial inversive distance circle packing metric to a metric with constant $\alpha$-curvature. More precisely, we prove that the inversive distance circle packing with constant $\alpha$-curvature is unique if $\alpha\chi(M)\leq 0$, which generalize Andreev-Thurston's rigidity results for circle packing with constant cone angles. We further prove that the solution to Ge-Xu's $\alpha$-flow can always be extended to a solution that exists for all time and converges exponentially fast to constant $\alpha$-curvature. Finally, we give some combinatorial and topological obstacles for the existence of constant $\alpha$-curvature metrics.
\end{abstract}


\section{Background}\label{Introduction}
This is a continuation of our study in \cite{Ge-Jiang1,Ge-Jiang2}. Koebe \cite{Koebe1} and Andreev \cite{Andreev2,Andreev1} first studied circle patterns (or say circle packings) with pairs of circles tangent to each other. Thurston generalized this to circle patterns with pairs of circles intersect with acute or right angles. However, Bowers and Stephenson \cite{Bowers-Stephenson} found it very useful to consider circle patterns with pairs of circles disjoint. In this case, they call it inversive distance circle patterns. In the pioneer work of \cite{CL1}, Chow and Luo first established an intrinsic connection between Thurston's circle pattern and surface Ricci flow. They use a discrete analog of Hamilton's Ricci flow to deform Thurston's circle pattern to one with constant cone angles. Inspired by Guo \cite{Guoren} and Luo's \cite{Luo1} extensive work on inversive distance circle patterns, and Bobenko, Pincall and Springborn's \cite{Bobenko} pioneer work on the extension of locally convex functional, we \cite{Ge-Jiang1,Ge-Jiang2} generalized Chow-Luo's discrete Ricci flow to inversive distance circle pattern setting. In \cite{Ge-Jiang1}, we consider inversive distance circle pattern in the Euclidean geometry background and deform it to one with constant cone angles. While in \cite{Ge-Jiang2}, we consider inversive distance circle pattern in the hyperbolic geometry background and deform it to one with flat cone angles. The classical definition of discrete Gaussian curvature (that is, $2\pi$ less cone angles) has intrinsic disadvantages. Ge and Xu \cite{Ge-Xu3} introduced a more reasonable discrete Gaussian curvature, which is more natural and respects the scaling exactly the same way as Gaussian curvature does. They \cite{Ge-Xu1} further introduced an $\alpha$-curvature, and an $\alpha$-flow which can be used to deform any initial circle pattern to the one with constant (or prescribed) $\alpha$-curvature. The $\alpha$-flow generalizes Chow-Luo's definition, in fact, the $\alpha=0$ case is exactly Chow-Luo's discrete Ricci flow. In this paper, we want to deform the inversive distance circle patterns along $\alpha$-flow to the one with constant $\alpha$-curvature. We also give some topological and combinatorial obstacles for the existence of constant $\alpha$-curvature metrics, which partially generalize Andreev and Thurston's results for circle patterns. We only consider the Euclidean background geometry in this paper.

Suppose $M$ is a closed surface with a triangulation $\mathcal{T}=\{V,E,F\}$,
where $V,E,F$ represent the sets of vertices, edges and faces respectively. Throughout this paper, all vertices are ordered one by one, marked by $1, \cdots, N$, where $N=|V|$ is the number of vertices. Throughout this paper, all functions $f: V\rightarrow \mathds{R}$ will be regarded as column vectors in $\mathds{R}^N$ and $f_i$ is the value of $f$ at $i$. We often write it as $f\in C(V)$. For fixed triangulated surface $(M, \mathcal{T})$, Thurston introduced circle packings (with acute or right intersection angles) to produce piecewise linear metrics on $M$, see \cite{T1} for more details. Bowers and Stephenson generalized Thurston's circle packing to inversive distance setting. For more explanations about inversive distance circle packings, see \cite{Bowers-Hurdal,Bowers-Stephenson,Guoren}. A inversive distance circle packing metric contains two factors. One is the inversive distance $I: E\rightarrow [0, +\infty)$ defined on all edges. The other is the circle radius vector $r: V\rightarrow (0, +\infty)$ defined on all vertices. For every given radius vector $r\in \mathds{R}^N_{>0}$ and inversive distance $I\geq0$, we equip each edge $\{ij\}\in E$ with an edge length
\begin{equation}\label{Def-edge-length}
l_{ij}=\sqrt{r_i^2+r_j^2+2r_ir_jI_{ij}}.
\end{equation}
However, for a face $\{ijk\}\in F$, the three lengths $l_{ij}, l_{jk}, l_{ik}$ may not satisfy the triangle inequalities any more. If for each $\{ijk\}\in F$, the three lengths $l_{ij}, l_{jk}, l_{ik}$ satisfy the triangle inequalities, then the triangulated surface $(M, \mathcal{T})$ with edge lengthes $l_{ij}$ forms an Euclidean polyhedral surface. Denote
\begin{equation}
\Omega\triangleq\Big\{r\in \mathds{R}^N_{>0}\;\big|\;l_{ij}+l_{jk}>l_{ik},\;l_{ij}+l_{ik}>l_{jk},\;l_{ik}+l_{jk}>l_{ij}, \;\forall \;\{ijk\}\in F\Big\}.
\end{equation}
In the following, each $r\in\Omega$ is called an inversive distance circle packing metric, or metric for short, and each $r\in\mathds{R}^N_{>0}$ is called a generalized inversive distance circle packing metric, or generalized metric for short. For each $r\in \Omega$, let $\theta_i^{jk}$ be the inner angle of triangle $\{ijk\}$
at the vertex $i$, the classical discrete Gauss curvature at the vertex $i$ is
\begin{equation}\label{Def-Gauss-curv}
K_i=2\pi-\sum_{\{ijk\}\in F}\theta_i^{jk},
\end{equation}
where the sum is taken over all the triangles with $i$ as one of its vertices.

\section{Definition of $\alpha$-curvature and $\alpha$-flow}\label{sec-defin-alpha-flow}
As is pointed by Ge-Xu \cite{Ge-Xu3}, there are two disadvantages of the classical definition of $K_i$. For one thing, classical discrete Gaussian curvature does not change under scaling, i.e. if $\tilde{r}_i=\lambda r_i$ for some positive constant $\lambda$, then $\tilde{K}_i=K_i$,
which is different from the transformation of scalar curvature $R_{\lambda g}=\lambda ^{-1}R_g$ in smooth case.
For another, classical discrete Gaussian curvature can't be used directly to approximate smooth Gaussian curvature since it always tends to zero as the triangulation becomes finer and finer. To amend this flaw, Ge-Xu suggested a new definition of discrete Gaussian curvature as $R_i=K_i/r_i^2$. $r_i^2$ seems to be the most reasonable discrete analog of smooth metric tensor $g$. Consider a smooth metric tensor $g$ with conical singularity at zero point, it can be expressed as $$g(z)=e^{f(z)}|z|^{2(\alpha-1)}dzd\bar{z}$$
locally for some $\alpha\in\mathds{R}$. Choosing $f(z)=-\ln\alpha^{2}$, then $g(z)=|dz^{\alpha}|^{2}$.
Comparing $r_i^{\alpha}$ with $|dz^{\alpha}|$, hence the $\alpha$-order metric $r_i^{\alpha}$ may be taken as a discrete analogue of conical metric to some extent. Based on this observation, Ge-Xu \cite{Ge-Xu3,Ge-Xu1} introduced the so called $\alpha$-curvature and $\alpha$-flow. We shall study the corresponding $\alpha$-curvature and $\alpha$-flow in the inversive distance circle packing setting.

\begin{definition}
Given $(M, \mathcal{T}, I)$ with inversive distance $I\geq0$. For any metric $r\in\Omega$ and any $\alpha\in\mathds{R}$, the discrete Gaussian curvature of order $\alpha$ (``$\alpha$-curvature" for short) at the vertex $i$ is defined to be
\begin{equation}\label{def-R-curvature}
R_{\alpha,i}=\frac{K_i}{r_i^{\alpha}},
\end{equation}
where $K_i$ is the classical discrete Gaussian curvature defined at $i$ by (\ref{Def-Gauss-curv}).
\end{definition}

In \cite{Ge-Xu1}, Ge and Xu introduced a $\alpha$-flow to deform any circle packing metrics to constant $\alpha$-curvature metrics. We consider the $\alpha$-flow in the inversive distance circle packing settings. Set $u_i=\ln r_i$.
\begin{definition}
Given $(M, \mathcal{T}, I)$ with inversive distance $I\geq0$. For any inversive distance circle packing metric $r\in\Omega$, denote $s_{\alpha}=2\pi\chi(M)/\|r\|^{\alpha}_{\alpha}$. The $\alpha$th order combinatorial Ricci (Yamabe) flow (``$\alpha$-flow" for short) is
\begin{equation}\label{Def-Yamabe-flow-2d}
\frac{du_i}{dt}=s_{\alpha}r_i^{\alpha}-K_i
\end{equation}
with initial value $u(0)\in\ln\Omega$, where $\ln\Omega$ is the image of $\Omega$ under the coordinate change $r\mapsto u$ with $u_i=\ln r_i$, $i\in V$.
\end{definition}

\begin{remark}
Chow and Luo's (normalized) combinatorial Ricci flow (see \cite{CL1}) is in fact the $\alpha$-flow defined above with $\alpha=0$ and $0\leq I\leq 1$. The flow studied in our former work \cite{Ge-Jiang1} is in fact the $\alpha$-flow with $\alpha=0$ and $I\geq0$.
\end{remark}

\section{Basic properties of $\alpha$-flow}\label{sec-basic-prop-alpha-flow}
We want to deform a inversive distance circle packing along the $\alpha$-flow. Since all terms in the $\alpha$-flow is smooth and hence locally Lipschitz continuous in $\Omega$, then by Picard theorem in classical ODE theory, the $\alpha$-flow (\ref{Def-Yamabe-flow-2d}) has a unique solution $u(t)$, $t\in[0, \epsilon)$ for some $\epsilon>0$.

\begin{proposition}\label{Prop-nonsingularity}
Given a triangulated surface $(M, \mathcal{T}, I)$ with inversive distance $I\geq0$. Let $u(t)$ be the solution to $\alpha$-flow (\ref{Def-Yamabe-flow-2d}), then the level set of $\prod_{i=1}^Nr_i(t)$ is preserved under the $\alpha$-flow. Furthermore, for each $i\in V$, $r_i(t)$ can not go to zero or infinity in any finite time interval $[0, a)$ with $a<+\infty$.
\end{proposition}
\noindent\begin{proof}
The proof is similar with \cite{Ge-Xu1}. Obviously, $(\sum u_i)'=\sum (s_{\alpha}r_i^{\alpha}-K_i)=0$. Then it follows that $\sum u_i(t)$ and hence $\prod r_i(t)$ are invariant along the $\alpha$-flow. Furthermore, note that $|s_{\alpha}r_i^{\alpha}|\leq2\pi|\chi(M)|$, $(2-d)\pi\leq K_i\leq 2\pi$, where $d=\max\limits_{1\leq i\leq N} d_i$ and $d_i$ is the degree of vertex $i$. Hence $|s_{\alpha}r_i^{\alpha}-K_i|$ is uniformly bounded by a constant $c>0$ depending only on the topology and the triangulation. Then
$$r_i(0)e^{-ct}\leq r_i(t)\leq r_i(0)e^{ct},$$
which implies that $r_i(t)$ can not go to zero or infinity in finite time.\qed
\end{proof}

\begin{remark}
We always assume the solution $r(t)$ satisfies $\prod_{i=1}^Nr_i(t)\equiv 1$ in the following. Denote $\mathscr{U}$ as the hyperplane $\{u\in \mathds{R}^N|\sum_{i=1}^Nu_i=0\}$, then $u(t)\subset\mathscr{U}$ as long as it exists.
\end{remark}

\begin{theorem}[$\alpha$-flow is variational]\label{Prop-gradient-flow}
The $\alpha$-flow (\ref{Def-Yamabe-flow-2d}) is a negative gradient flow.
\end{theorem}
\noindent\begin{proof}
Consider a generic triangle $\{ijk\}\in F$, which is configured by three circles with nonnegative inversive distance $I_{ij}$, $I_{jk}$ and $I_{ik}$. Recall the notation below
\begin{equation}
\Delta_{ijk}\triangleq\left\{(r_i,r_j,r_k)\in\mathds{R}^3_{>0}\,\big|\,l_{ij}+l_{jk}>l_{ik},l_{ij}+l_{ik}>l_{jk},l_{jk}+l_{ik}>l_{ij}\right\}.
\end{equation}
For each metric $r\in\Omega$, it is obviously that $(r_i,r_j,r_k)\in \Delta_{ijk}$ for any triangle $\{ijk\}\in F$. We use $\ln \Delta_{ijk}$ as the homeomorphic image of $\Delta_{ijk}$ under coordinate change $u_i=\ln r_i$. Guo \cite{Guoren} once proved that $\partial \theta_i^{jk}/\partial u_j=\partial \theta_j^{ik}/\partial u_i$ in $\ln \Delta_{ijk}$ and $\ln \Delta_{ijk}$ is a simply connected nonvoid open set. Hence $\theta_i^{jk}du_i+\theta_j^{ik}du_j+\theta_k^{ij}du_k$ is a closed $C^{\infty}$-smooth $1$-form in $\ln \Delta_{ijk}$. For any fixed $u_0=(u_{0,1},\cdots,u_{0,N})^T\in\ln \Omega$, then for each $(u_i,u_j,u_k)\in\ln \Delta_{ijk}$, the following line integration
\begin{equation}
F_{ijk}(u_i,u_j,u_k)\triangleq\int_{(u_{0,i},u_{0,j},u_{0,k})}^{(u_i,u_j,u_k)}\theta_i^{jk}du_i+\theta_j^{ik}du_j+\theta_k^{ij}du_k
\end{equation}
is well defined. Set
\begin{equation}
B_{\alpha}(u)=\frac{2\pi\chi(M)}{\alpha}\Big(\ln\big(e^{\alpha u_1}+\cdots+e^{\alpha u_N}\big)\Big)\Big|^u_{u_0}.
\end{equation}
It's easy to see $\partial B_{\alpha}(u)/\partial u_i=s_{\alpha}r_i^{\alpha}$, moreover, $B_{\alpha}(u)$ depends on $\alpha$ continuously, that is
\begin{equation*}
\lim\limits_{\alpha\rightarrow 0}B_{\alpha}(u)=B_0(u).
\end{equation*}
Therefore we can define the $\alpha$-order discrete Ricci potential (``$\alpha$-potential" for short) as
\begin{equation}\label{def-alpha-potential}
F(u)=A(u)-B_{\alpha}(u)-\sum_{\{ijk\}\in F}F_{ijk}(u_i,u_j,u_k),
\end{equation}
where $A(u)=2\pi\sum_{i=1}^N(u_i-u_{0,i})$ is a linear function of $u$. By calculation, we get
\begin{equation}
\frac{\partial F}{\partial u_i}=K_i-s_{\alpha}r_i^{\alpha}.
\end{equation}
Thus the $\alpha$-flow is in fact the negative gradient flow $\dot{u}=-\nabla_u F$ of the $\alpha$-potential $F$.\qed
\end{proof}

\begin{lemma}\label{L-semi-positive-definite}
Given a triangulated surface $(M, \mathcal{T})$ with inversive distance $I\geq 0$. Denote
\begin{equation}
L\triangleq\frac{\partial(K_{1},\cdots,K_{N})}{\partial(u_{1},\cdots,u_{N})}.
\end{equation}
Then $L$ is positive semi-definite with rank $N-1$. Moreover, the kernel of $L$ is the linear space spanned by the vector $\mathds{1}\triangleq(1,\cdots,1)^T$.
\end{lemma}
\begin{remark}
For $I\in[0,1]$ case, which is equivalent to say Andreev-Thurston's circle pattern with weight $\Phi\in[0,\frac{\pi}{2}]$ case, Chow-Luo proved above lemma, see Proposition 3.9 in \cite{CL1}. For the general inversive distance $I\geq0$ case, Guo proved above lemma essentially. See Lemma 5, Lemma 6 and Corollary 7 in \cite{Guoren}, or see Proposition 3.8 in \cite{Ge-Jiang1} for a direct proof.
\end{remark}
Recall Ge-Xu \cite{Ge-Xu1} once introduced an $\alpha$-order combinatorial Lalacian
\begin{equation}\label{Def-alpha-Laplacian-2d}
\Delta_\alpha f_i=-\frac{1}{r_i^\alpha}\sum_{j=1}^N \frac{\partial K_i}{\partial u_j} f_j=-\frac{1}{r_i^\alpha}\sum_{j\thicksim i} \frac{\partial K_i}{\partial u_j}(f_j-f_i)
\end{equation}
which is still meaningful in the inversive distance circle packing setting. Denote
$$r^{\alpha}=(r_1^{\alpha},\cdots,r_N^{\alpha})^T,$$
and $\Sigma=diag\big\{r_1,\cdots,r_N\big\}$. Note each $f\in C(V)$ is also a column vector, $\alpha$-Laplacian (\ref{Def-alpha-Laplacian-2d}) can be written in a matrix form,
\begin{equation}
\Delta_{\alpha}=-\Sigma^{-\alpha}L
\end{equation}
with $\Delta_{\alpha} f=-\Sigma^{-\alpha}Lf$ for each $f\in C(V)$. $\alpha$-order Combinatorial Laplace operator has similar properties with smooth Laplace operator, for details see \cite{Ge-Xu1,Ge-Xu2}. Denote the first positive eigenvalue of $-\Delta_{\alpha}$ by $\lambda_1(-\Delta_{\alpha})$.

\begin{theorem}\label{Thm-converg-exist-const-cur-metric}
Given a triangulated surface $(M, \mathcal{T})$ with inversive distance $I\geq 0$. If the solution $\{r(t)\}\subset\Omega$ to $\alpha$-flow (\ref{Def-Yamabe-flow-2d}) converges to a metric $r^*\in\Omega$, then $r^*$ is a constant $\alpha$-curvature metric. As a consequence, there is a constant $\alpha$-curvature metric on $(M, \mathcal{T})$.
\end{theorem}
\noindent\begin{proof}
We say that the solution $r(t)$ to $\alpha$-flow (\ref{Def-Yamabe-flow-2d}) converges, if $r(t)$ exists for all time $t\in[0,+\infty)$, and $r(t)$ converges in the Euclidean space topology to some $r^*\in\Omega$ as time $t$ tends to $+\infty$. According to the classical ordinary differential equation theory, for each vertex $i\in V$, $r^*$ is a zero point of $(s_{\alpha}r_i^{\alpha}-K_i)r_i$. Hence $K_i/r_i^{\alpha}\equiv s_{\alpha}$ for each $i$, implying that $r^*$ is a constant
$\alpha$-curvature metric.\qed\\
\end{proof}

Conversely, we have
\begin{theorem}[Local stable]\label{Thm-const-curv-metric-exist-converg}
Given a triangulated surface $(M, \mathcal{T})$ with inversive distance $I\geq 0$. Let $r^*\in\Omega$ be a constant $\alpha$-curvature metric. Assume $\lambda_1(-\Delta_{\alpha})>\alpha s_{\alpha}$ at $r^*$, then the solution $r(t)$ to $\alpha$-flow (\ref{Def-Yamabe-flow-2d}) exists for all time $t\geq0$ and converges exponentially fast to $r^*$ if $r(0)$ is close enough to $r^*$.
\end{theorem}
\noindent\begin{proof}
Let $u^*$ be the corresponding $u$-coordinate of $r^*$. We follow the proof of Theorem 1 in \cite{Ge-Xu1}. Write the $\alpha$-flow (\ref{Def-Yamabe-flow-2d}) in a matrix form as
$$\dot{u}=s_{\alpha}r^{\alpha}-K,$$
then differentiate $s_{\alpha}r^{\alpha}-K$ at $u^*$, we get
\begin{equation}
-D_u(s_{\alpha}r^{\alpha}-K)\big|_{u^*}=L-\alpha s_{\alpha}
\left(\Sigma^{\alpha}-\frac{r^{\alpha}(r^{\alpha})^T}{\|r\|_{\alpha}^{\alpha}}\right).
\end{equation}
Set $\Lambda_{\alpha}=\Sigma^{-\frac{\alpha}{2}}L\Sigma^{-\frac{\alpha}{2}}$, then
$$\Delta_{\alpha}=-\Sigma^{-\frac{\alpha}{2}}\Lambda_{\alpha}\Sigma^{\frac{\alpha}{2}},$$
which implies that $\lambda_1(-\Delta_{\alpha})=\lambda_1(\Lambda_{\alpha})$. For $u\in\ln\Omega$, it's easy to calculate
\begin{equation}\label{Hession of F}
\begin{aligned}
L-\alpha s_{\alpha}
\left(\Sigma^{\alpha}-\frac{r^{\alpha}(r^{\alpha})^T}{\|r\|_{\alpha}^{\alpha}}\right)=\Sigma^{\frac{\alpha}{2}}
\left(\Lambda_{\alpha}-\alpha s_{\alpha}\left(I-\frac{r^{\frac{\alpha}{2}}(r^{\frac{\alpha}{2}})^T}{\|r\|^{\alpha}_{\alpha}}\right)\right)
\Sigma^{\frac{\alpha}{2}}.
\end{aligned}
\end{equation}
Choose an orthonormal matrix $P$ such that $P^T\Lambda_{\alpha} P=diag\{0,\lambda_1(\Lambda_{\alpha}),\cdots,\lambda_{N-1}(\Lambda_{\alpha})\}$.
Suppose $P=(e_0,e_1,\cdots,e_{N-1})$,
where $e_i$ is the $(i+1)$-column of $P$.
Then $\Lambda_{\alpha} e_0=0$ and $\Lambda_{\alpha} e_i=\lambda_i e_i,\,1\leq i\leq N-1$,
which implies $e_0=r^\frac{\alpha}{2}/\|r^\frac{\alpha}{2}\|$ and $r^\frac{\alpha}{2}\perp e_i,\,1\leq i\leq N-1$.
Hence $\big(I-\frac{r^{\frac{\alpha}{2}}(r^{\frac{\alpha}{2}})^T}{\|r\|^{\alpha}_{\alpha}}\big)e_0=0$ and $\big(I-\frac{r^{\frac{\alpha}{2}}(r^{\frac{\alpha}{2}})^T}{\|r\|^{\alpha}_{\alpha}}\big)e_i=e_i$, $1\leq i\leq N-1$,
thus we get
\begin{equation*}
\begin{aligned}
-D_u(s_{\alpha}r^{\alpha}-K)\big|_{u^*}=\Sigma^{\frac{\alpha}{2}}P
\left(
       \begin{array}{cccc}
       0 & & &   \\
       &  \lambda_1(\Delta_{\alpha})-\alpha s_{\alpha} &   &0  \\
       & & \ddots &   \\
       & 0& & \lambda_{N-1}(\Delta_{\alpha})-\alpha s_{\alpha} \\
       \end{array}
\right)P^T\Sigma^{\frac{\alpha}{2}}.
\end{aligned}
\end{equation*}
If $\lambda_1(-\Delta_{\alpha})>\alpha s_{\alpha}$ at $r^*$, then $D_u(s_{\alpha}r^{\alpha}-K)\big|_{u^*}$ is semi-negative definite with rank $N-1$. Moreover, its kernel is spanned by vector $\mathds{1}$, which is perpendicular to the hyperplane $\mathscr{U}$. Hence $D_u(s_{\alpha}r^{\alpha}-K)\big|_{u^*}$ is strictly negative definite when restricted in the hyperplane $\mathscr{U}$. By Lyapunov stability theorem in ODE theory, $u^*$
is an asymptotically stable point of the $\alpha$-flow (\ref{Def-Yamabe-flow-2d}), this implies the conclusion above.\qed\\
\end{proof}

If $\alpha\chi(M)\leq 0$, then there always $\lambda_1(-\Delta_{\alpha})>0>\alpha s_{\alpha}$. Therefore we get
\begin{corollary}\label{corollary-const-curv-metric-exist-converg}
Given $(M, \mathcal{T})$ with inversive distance $I\geq 0$. Let $r^*\in\Omega$ be a constant $\alpha$-curvature metric. If $\alpha\chi(M)\leq 0$, then the solution $r(t)$ to $\alpha$-flow (\ref{Def-Yamabe-flow-2d}) exists for all time $t\geq0$ and converges exponentially fast to $r^*$ if $r(0)$ is close enough to $r^*$.
\end{corollary}

\section{Extended $\alpha$-flow}\label{sec-extend-alpha-flow}
We follow the idea from \cite{Bobenko,Luo1,Luo2} to give the definition of generalized triangle. A generalized Euclidean triangle $\triangle 123$ is a (topological) triangle of vertices $v_1, v_2, v_3$ so that each edge is assigned a positive number, called edge length. Let $x_i$ be the assigned length of the edge $v_jv_k$ where $\{i, j, k\}$=$\{1, 2, 3\}$. The inner angle $\tilde{\theta}_i$=$\tilde{\theta}_i(x_1, x_2, x_3)$ at the vertex $v_i$ is defined as follows. If $x_1, x_2, x_3$ satisfy the triangle inequalities that $x_j+x_k>x_i$ for $\{i, j, k\}$=$\{1, 2, 3\}$, then $\tilde{\theta}_i$ is the inner angle of the corresponding Euclidean triangle of edge lengths $x_1, x_2, x_3$ opposite to the edge of length $x_i$; if $x_i\geq x_j+x_k$, then $\tilde{\theta}_i=\pi$, and $\tilde{\theta}_j=\tilde{\theta}_k=0$. For each $r\in\mathds{R}^N_{>0}$ and each $\{ijk\}\in F$, $l_{ij}$, $l_{jk}$ and $l_{ik}$ form the three edge lengthes of a generalized Euclidean triangle. Hence the inner angle map $\tilde{\theta}_i^{jk}=\tilde{\theta}_i^{jk}(r)$ is defined for any $r\in\mathds{R}^N_{>0}$. Luo \cite{Luo1} proved that it is a continuous extension of function $\theta_i^{jk}=\theta_i^{jk}(r)$, $r\in\Omega$. Using the auxiliary function
\begin{equation}
\Lambda(x)=
\begin{cases}
\;\;\;\;\;\pi, & \text{$x \leq -1$}\\
\arccos x, & \text{$-1\leq x \leq 1$}\\
\;\;\;\;\;0,& \text{$x\geq 1$}
\end{cases}
\end{equation}
we introduced in \cite{Ge-Jiang1}, we can express the generalized inner angle $\tilde{\theta}_i^{jk}$ as
\begin{equation}\label{formula-tilde-xita-euclid}
\tilde{\theta}_i^{jk}(r_i,r_j,r_k)=\Lambda\bigg(\frac{l_{ij}^2+l_{ik}^2-l_{jk}^2}{2l_{ij}l_{ik}}\bigg),
\end{equation}
which also implies that $\tilde{\theta}_i^{jk}$ is a continuous extension of $\theta_i^{jk}$. It follows that the definition of $K_i$, which is originally defined on $r\in\Omega$, now can be  naturally extended to
\begin{equation}\label{def-K-tuta-alpha-flow}
\widetilde{K}_i=2\pi-\sum_{\{ijk\} \in F}\tilde{\theta}_i^{jk}.
\end{equation}
Hence the curvature map $K(r):\Omega\rightarrow \mathds{R}^N_{>0}$ is extended continuously to $\widetilde{K}(r):\mathds{R}^N_{>0}\rightarrow \mathds{R}^N$
with discrete Gauss-Bonnet formula
\begin{equation}\label{Gauss-Bonnet-extend}
\sum_{i=1}^N\widetilde{K}_i=2\pi \chi(M).
\end{equation}
valid, see Proposition 3.5 in \cite{Ge-Jiang1}.

\begin{definition}[extended $\alpha$-flow]
Given a triangulated surface $(M, \mathcal{T}, I)$ with inversive distance $I\geq 0$. The extended $\alpha$-flow is defined as for each $i\in V$,
\begin{equation}\label{Def-extended-flow-u}
\frac{du_i}{dt}=s_{\alpha}r_i^{\alpha}-\widetilde{K}_i.
\end{equation}
\end{definition}

\begin{theorem}\label{Thm-main-1}
Given a triangulated surface $(M, \mathcal{T}, I)$ with inversive distance $I\geq 0$. Suppose $\{u(t)|t\in[0,T)\}$ is the unique maximal solution to $\alpha$-flow (\ref{Def-Yamabe-flow-2d}) with $0<T\leq +\infty$. Then we can always extend it to a solution $\{u(t)|t\in[0,+\infty)\}$ when $T<+\infty$, that is, the solution to the extended $\alpha$-flow (\ref{Def-extended-flow-u}) exists for all time $t\in[0,+\infty)$ for any $u(0)\in \ln\Omega$.
\end{theorem}
\noindent\begin{proof}
The proof is similar with Proposition \ref{Prop-nonsingularity}. Obviously, the extended curvature $\widetilde{K}$ satisfies $(2-d)\pi\leq K_i\leq 2\pi$ too. Hence $|s_{\alpha}r_i^{\alpha}-\widetilde{K}_i|$ is uniformly bounded too, this implies the conclusion above by the extension theorem for solution in ODE theory.\qed
\end{proof}
\begin{remark}
We don't know whether the extended solution $\{u(t)|t\in[0,+\infty)\}$ is unique when $T<+\infty$. Similar to what is shown in subsection 5.2 of \cite{Ge-Jiang2}, in the Euclidean geometry background, when ever $u(t)$ lies in $\ln\Omega$, it is unique since all the coefficients are smooth and then locally Lipschitz continuous in $\ln\Omega$. However, $\widetilde{K}_i$ is generally not Lipschitz continuous near the ``triangle inequality invalid boundary" of $\ln\Omega$. This is because
$\partial \theta_i^{jk}/\partial r_j=\infty$ at the ``triangle inequality invalid boundary". So we don't know whether the extended solution is unique.
\end{remark}

\section{Uniqueness of constant $\alpha$-curvature metric}\label{sec-uniqe-const-alpha-cur-metric}
\begin{theorem}\label{thm-uniqueness-for-alpha-curvature}
Given a triangulated surface $(M, \mathcal{T})$ with inversive distance $I\geq 0$. Assume $\alpha\chi(M)\leq 0$, then the constant $\alpha$-curvature metric is unique (up to scaling) if it exists.
\end{theorem}

\begin{remark}
Before giving the proof of Theorem \ref{thm-uniqueness-for-alpha-curvature}, we want to remark, there are evidences (such as Example 3 in \cite{Ge-Xu3}) showing that the constant $\alpha$-curvature metrics may be not unique if not assuming $\alpha\chi(M)\leq 0$. And in this sense, the result in Theorem \ref{thm-uniqueness-for-alpha-curvature} is sharp.
\end{remark}

\noindent\begin{proof}
We shall prove, if there is a metric (in $u$-coordinate) $p\in\ln \Omega\cap\mathscr{U}$ such that the corresponding $\alpha$-curvature is constant, then all other (generalized) metrics in $\mathscr{U}$ can't have constant $\alpha$-curvature any more. We first prove the $C^1$-smooth and convex extension property of $\alpha$-potential (\ref{def-alpha-potential}).
\begin{claim}\label{claim-hessian-F}
If $\alpha\chi(M)\leq 0$, then the $\alpha$-potential (\ref{def-alpha-potential}) can be extended to a convex function $\widetilde{F}$ defined on $\mathds{R}^N$ with $\widetilde{F}\in C^1(\mathds{R}^N)\cap C^{\infty}(\ln \Omega)$ and $\widetilde{F}(u)=\widetilde{F}(u+t\mathds{1})$ for any $t\in\mathds{R}$. Moreover, in $\ln\Omega$, $Hess_u\widetilde{F}$ is positive semi-definite with rank $N-1$ and kernel $\{t\mathds{1}|t\in\mathds{R}\}$.
\end{claim}
By Luo's extension theory (see Proposition 3.2 in \cite{Luo1}, or Lemma 3.6 in \cite{Ge-Jiang1}), for each triangle $\{ijk\}\in F$, and for any $(u_i,u_j,u_k)\in\mathds{R}^3$, the following integration
\begin{equation}
\widetilde{F}_{ijk}(u_i,u_j,u_k)\triangleq\int_{(u_{0,i},u_{0,j},u_{0,k})}^{(u_i,u_j,u_k)}\tilde{\theta}_i^{jk} du_i+\tilde{\theta}_j^{ik} du_j+\tilde{\theta}_k^{ij}du_k
\end{equation}
is well defined and is a $C^1$-smooth concave function on $\mathds{R}^3$. Moreover, $\widetilde{F}_{ijk}(u_i,u_j,u_k)$ is the $C^1$-extension of $F_{ijk}(u_i,u_j,u_k)$, from $(u_i,u_j,u_k)\in \ln\Delta_{ijk}$ to $(u_i,u_j,u_k)\in\mathds{R}^3$. Then the extended $\alpha$-potential $\widetilde{F}(u)$ can be defined as
\begin{equation}\label{def-extended-alpha-potential}
\widetilde{F}(u)=A(u)-B_{\alpha}(u)-\sum_{\{ijk\}\in F}\widetilde{F}_{ijk}(u),
\end{equation}
It's easy to see, if $\alpha\chi(M)\leq 0$, then $B_{\alpha}(u)$ is convex on the whole space $\mathds{R}^N$ and is strictly convex when constrained to the hyperplane $\mathscr{U}$. Hence $\widetilde{F}(u)$ is convex extension of $F$ and is $C^1$-smooth in $\mathds{R}^N$. Using discrete Gauss-Bonnet formula (\ref{Gauss-Bonnet-extend}), one see $\sum_{i=1}^N(\widetilde{K}_i-s_{\alpha}r_i^{\alpha})\equiv 0$, this implies $\widetilde{F}(u)=\widetilde{F}(u+t\mathds{1})$ for any $t\in\mathds{R}$ and any $u\in\mathds{R}^N$. In $\ln \Omega$, $\widetilde{F}=F$ and hence is $C^{\infty}$-smooth.  then it's easy to calculate
\begin{equation}\label{Hession of F}
\begin{aligned}
Hess_u\widetilde{F}=Hess_uF
=L-\alpha s_{\alpha}
\left(\Sigma^{\alpha}-\frac{r^{\alpha}(r^{\alpha})^T}{\|r\|_{\alpha}^{\alpha}}\right),
\end{aligned}
\end{equation}
on $\ln\Omega$, where $L=\frac{\partial(K_1, \cdots, K_N)}{\partial(u_1, \cdots, u_N)}$. Obviously, both matrix $\Sigma^{\alpha}-\frac{r^{\alpha}(r^{\alpha})^T}{\|r\|_{\alpha}^{\alpha}}$ and $L$ (see Lemma \ref{L-semi-positive-definite})
are positive semi-definite with rank $N-1$ and null space $\{t\mathds{1}|t\in\mathds{R}\}$, which implies Claim \ref{claim-hessian-F}. Next we prove
\begin{claim}\label{claim-F-2}
Assuming $\alpha\chi(M)\leq 0$ and there is a metric (in $u$-coordinate) $p\in\ln \Omega\cap\mathscr{U}$ such that the corresponding $\alpha$-curvature is constant. For each direction $\xi\in \mathbb{S}^{N-1}\cap\mathscr{U}$, set $\varphi_{\xi}(t)=\widetilde{F}(p+t\xi)$, $t\in \mathds{R}$. Then $\varphi_{\xi}'(0)=0$. Moreover, $\varphi_{\xi}'(t)$ is increasing on $\mathds{R}$ and is strictly increasing on $[-\delta, \delta]$, for any $\delta$ with $0<\delta<dist(p,\partial\ln\Omega)$.
\end{claim}
In fact, $\nabla_u\widetilde{F}=\widetilde{K}-s_{\alpha}r^{\alpha}=(\widetilde{K}_1-s_{\alpha}r^{\alpha}_1,\cdots,\widetilde{K}_N-s_{\alpha}r^{\alpha}_N)^T$.
Hence $\varphi_{\xi}'(t)=(\widetilde{K}-s_{\alpha}r^{\alpha})\cdot\xi$. Note for $t\in[-\delta, \delta]$, $p+t\xi\in\ln\Omega\cap\mathscr{U}$, hence
$\varphi_{\xi}'(t)=(K-s_{\alpha}r^{\alpha})\cdot\xi$ and then $\varphi_{\xi}'(0)=0$. On one hand, $\widetilde{F}$ is convex on $\mathds{R}^N$, then $\varphi_{\xi}(t)$ is convex on $\mathds{R}$, which implies $\varphi_{\xi}'(t)$ is increasing for $t\in\mathds{R}$. On the other hand, by Claim \ref{claim-hessian-F}, for $u\in B(p,\delta)$, the kernel space of $Hess_uF$ is exactly perpendicular to the hyperplane $\mathscr{U}$, hence $F$ is strictly convex on $B(p,\delta)\cap\mathscr{U}$, then $\varphi_{\xi}(t)$ is strictly convex on $(-\delta,\delta)$, which implies that  $\varphi_{\xi}'(t)$ is strictly increasing on $(-\delta,\delta)$. Thus we get Claim \ref{claim-F-2}.

Claim \ref{claim-F-2} implies that $t=0$ is the unique zero point of $\varphi_{\xi}'(t)$. Using this fact, we can finally prove the theorem.
Under the assumption of Claim \ref{claim-F-2}, there are no other point in $\mathscr{U}$ with constant $\alpha$-curvature,
because $\varphi_{\xi}'=0$ at each constant $\alpha$-curvature metric (may be a generalized metric). \qed

\end{proof}

\section{Deform the $\alpha$-curvature to a constant}\label{sec-defom-alpha-cur-to-const}
\begin{lemma}\label{lemma-gotoinfinity}
Assuming $f\in C(\mathds{R}^n)$ and for any direction $\xi\in \mathbb{S}^{n-1}$, $f(t\xi)$ as a function of $t$ is monotone increasing on $[0,+\infty)$ and tends to $+\infty$ as $t\rightarrow +\infty$. Then $\lim\limits_{x\rightarrow\infty}f(x)=+\infty.$
\end{lemma}
For a proof of above lemma, see \cite{Ge-Jiang1}. Using this lemma, we finally get

\begin{theorem}[Long term convergence]\label{Thm-main-2}
Given $(M, \mathcal{T})$ with inversive distance $I\geq 0$. Assume $\alpha\chi(M)\leq 0$ and there is a metric of constant $\alpha$-curvature $r^*\in\Omega$. Then $r(t)$ can always be extended to a solution that converges exponentially fast to a constant $\alpha$-curvature metric as $t\rightarrow+\infty$. That is, any solution to the extended flow
(\ref{Def-extended-flow-u}) converges exponentially fast to a metric of constant $\alpha$-curvature as $t\rightarrow+\infty$.
\end{theorem}
\noindent\begin{proof}
By Claim \ref{claim-F-2}, $\varphi'_{\xi}(t)\geq\varphi'_{\xi}(\delta)>0$ for $t>\delta$ while $\varphi'_{\xi}(t)\leq\varphi'_{\xi}(-\delta)<0$ for $t<-\delta$. Hence
$$\varphi_{\xi}(t)\geq\varphi_{\xi}(\delta)+\varphi'_{\xi}(\delta)(t-\delta)$$
for $t\geq \delta$, while
$$\varphi_{\xi}(t)\geq\varphi_{\xi}(-\delta)+\varphi'_{\xi}(-\delta)(t+\delta)$$
for $t\leq -\delta$. This implies
\begin{equation}
\lim\limits_{t\rightarrow\pm\infty}\varphi_{\xi}(t)=+\infty.
\end{equation}
Using Lemma \ref{lemma-gotoinfinity}, we know that $\widetilde{F}(u)$ is proper on $\mathscr{U}$ and
$\lim\limits_{u\in \mathscr{U}, u\rightarrow\infty}\widetilde{F}(u)=+\infty.$

Let $\varphi(t)=\widetilde{F}(u(t))$, then $\varphi'(t)=-\|\widetilde{K}-s_{\alpha}r^{\alpha}\|^2\leq 0$. Thus $\varphi(t)$ is decreasing, note $\widetilde{F}$ is proper, then $\varphi(t)$ is compactly supported in $\mathds{R}^N$. Also, $\widetilde{F}$ is bounded form below, hence $\varphi(+\infty)$ exists. By the mean value theorem, there exists a sequence $t_n\in(n,n+1)$ such that $\varphi'(t_n)=\varphi(n+1)-\varphi(n)\rightarrow 0$. By choosing a subsequence of $t_n$, which is still denoted by $t_n$, we require that $u(t_n)$ converges to some point $\hat{u}$. Hence $\widetilde{K}(u(t_n))$ converges to $\widetilde{K}(\hat{u})$. Using $\varphi'(t_n)=-\|\widetilde{K}(u(t_n))-s_{\alpha}r^{\alpha}\|^2\rightarrow 0$, we have $\widetilde{K}(\hat{u})=s_{\alpha}\hat{r}^{\alpha}$. By the uniqueness of constant $\alpha$-curvature metric proved in Theorem \ref{thm-uniqueness-for-alpha-curvature}, we get $u^*=\hat{u}$ and $u(t_n)\rightarrow u^*$. It means for some sufficient big $t_{n}$, $u(t_{n})$ is very close to constant curvature metric $u^*$. Further note that in Corollary \ref{corollary-const-curv-metric-exist-converg}, we had already proved that $u^*$ is the asymptotically stable point of the $\alpha$-flow. Then the solution $\{u(t)\}_{t\geq t_n}$ converges exponentially fast to $u^*$, i.e., the original solution $\{u(t)\}_{t\geq 0}$ converges exponentially fast to $u^*$.\qed\\
\end{proof}

\begin{remark}
It's easy to see, the $\alpha$-potential (\ref{def-alpha-potential}) can always be extended to a function $\widetilde{F}$ defined on $\mathds{R}^N$ with $\widetilde{F}\in C^1(\mathds{R}^N)\cap C^{\infty}(\ln \Omega)$ and $\widetilde{F}(u)=\widetilde{F}(u+t\mathds{1})$ for any $t\in\mathds{R}$. If further assume $\lambda_1(\Delta_{\alpha})>\alpha s_{\alpha}$ in $\ln\Omega$, then $Hess_uF\geq0$, $rank(Hess_uF)=N-1$, the null space of $Hess_uF$ is $\{t\mathds{1}|t\in\mathds{R}\}$. This implies that $\widetilde{F}$ is locally convex in $\ln\Omega$ and is locally strictly convex in $\ln\Omega\cap\mathscr{U}$. Generally, $\widetilde{F}$ is not convex in $\mathds{R}^N$, see Example 2 in \cite{Ge-Xu3}. While Claim \ref{claim-hessian-F} says $\widetilde{F}$ is indeed convex in $\mathds{R}^N$ when assuming $\alpha\chi(M)\leq 0$. Obviously, condition $\alpha\chi(M)\leq 0$ is much stronger than and hence implies condition $\lambda_1(-\Delta_{\alpha})>\alpha s_{\alpha}$ in $\ln\Omega$. We expect $\widetilde{F}$ still convex in $\mathds{R}^N$ under the assumption $\lambda_1(\Delta_{\alpha})>\alpha s_{\alpha}$ in $\ln\Omega$. If so, then Theorem \ref{thm-uniqueness-for-alpha-curvature} and Theorem \ref{Thm-main-2} are still valid if we substitute condition $\lambda_1(-\Delta_{\alpha})>\alpha s_{\alpha}$ everywhere in $\ln\Omega$ for condition $\alpha\chi(M)\leq 0$.
\end{remark}

\section{Combinatorial-Topological obstacle}
We want to know the interplay between combinatorial-topological structure of $(M,\mathcal{T})$ and the constant $\alpha$-curvature metric. For any nonempty proper subset $A\subset V$, let $F_A$ be the subcomplex whose vertices are in $A$ and let $Lk(A)$ be the set of pairs $(e, v)$ of an edge $e$ and a vertex $v$ satisfying the following three conditioins: (1) The end points of $e$ are not in $A$; (2) $v$ is in $A$; (3) $e$ and $v$ form a triangle. For any nonempty proper subset $A\subset V$, the notation
\begin{equation*}
Y_A\triangleq\Big\{x\in \mathds{R}^N \Big|\sum_{i\in A}x_i >-\sum_{(e,v)\in Lk(A)}\big(\pi-\Lambda(I_e)\big)+2\pi\chi(F_A)\Big\},
\end{equation*}
defined in \cite{Ge-Jiang1}. Write
\begin{equation}\label{convex-set}
Y\triangleq\Big\{x\in \mathds{R}^N\Big|\sum_{i\in V}x_i=2\pi\chi(M)\Big\}\bigcap\Big(\mathop{\bigcap} _{\phi\neq A\subsetneqq V} Y_A \Big).
\end{equation}
We combine Theorem 5.4 and Theorem 5.8 in \cite{Ge-Jiang1} together and restate it here as
\begin{theorem}\label{thm-image-domain-K}
Given a triangulated surface $(M, \mathcal{T})$ with inversive distance $I\geq 0$. Then the space of all possible discrete Gaussian curvatures $K(\Omega)$ is contained in $Y$, while the space of all possible extended curvatures $\widetilde{K}(\mathds{R}^N_{>0})$ is contained in $\overline{Y}$, the closure of $Y$.
\end{theorem}

Using Theorem \ref{thm-image-domain-K}, we get the following four corollaries directly.
\begin{corollary}\label{cor-combin-topo-1}
Given a triangulated surface $(M, \mathcal{T}, I)$ with inversive distance $I\geq 0$. If $r^*\in\Omega$ is a metric with constant $\alpha$-curvature, then for any nonempty proper subset $A\subset V$,
\begin{equation*}
2\pi\chi(M)\frac{\sum_{i\in A}r_i^*}{\|r^*\|_{\alpha}^{\alpha}} >-\sum_{(e,v)\in Lk(A)}\big(\pi-\Lambda(I_e)\big)+2\pi\chi(F_A).
\end{equation*}
\end{corollary}

\begin{corollary}\label{cor-combin-topo-2}
Given a triangulated surface $(M, \mathcal{T}, I)$ with inversive distance $I\geq 0$. If $r^*\in\mathds{R}_{>0}^N$ is a (extended) metric with constant (extended) $\alpha$-curvature, then for any nonempty proper subset $A\subset V$,
\begin{equation*}
2\pi\chi(M)\frac{\sum_{i\in A}r_i^*}{\|r^*\|_{\alpha}^{\alpha}} \geq-\sum_{(e,v)\in Lk(A)}\big(\pi-\Lambda(I_e)\big)+2\pi\chi(F_A).
\end{equation*}
\end{corollary}

\begin{corollary}\label{cor-combin-topo-3}
Given $(M, \mathcal{T}, I)$ with inversive distance $I\geq 0$. Assuming there exists a metric in $\Omega$ with constant $\alpha$-curvature. If $\chi(M)>0$, then
$Y\cap\mathds{R}_{>0}^N\neq\emptyset$; If $\chi(M)<0$, then $Y\cap\mathds{R}_{<0}^N\neq\emptyset$; If $\chi(M)=0$, then $Y\cap\{0\}\neq\emptyset$.
\end{corollary}

\begin{corollary}\label{cor-combin-topo-4}
Given $(M, \mathcal{T}, I)$ with inversive distance $I\geq 0$. Assuming there exists a extended metric in $\mathds{R}_{>0}^N$ with constant extended $\alpha$-curvature. If $\chi(M)>0$, then $\overline{Y}\cap\mathds{R}_{>0}^N\neq\emptyset$; If $\chi(M)<0$, then $\overline{Y}\cap\mathds{R}_{<0}^N\neq\emptyset$; If $\chi(M)=0$, then $\overline{Y}\cap\{0\}\neq\emptyset$.
\end{corollary}

Obviously, $Y$ is exactly determined by the topological information of surface $M$ and combinatorial information of triangulation $\mathcal{T}$. Thus Corollary \ref{cor-combin-topo-1}-\ref{cor-combin-topo-4} encode the existence of constant (extended) $\alpha$-curvature metrics with combinatorial-topological obstacles.
As was proved by the first author and Xu \cite{Ge-Xu3}, in Thurston's circle packing setting, which is equivalent to say $0\leq I\leq1$, these combinatorial-topological obstacles are the only obstacles for the existence of constant $\alpha$-curvature metric when assuming $\alpha\chi(M)\leq0$ (see Theorem 1.5, Theorem 1.6 in \cite{Ge-Xu3}). In inversive distance circle packing setting, things may be different.  \\

\noindent \textbf{Acknowledgements}: We would like to thank BICMR and Professor Gang Tian for constant support and encouragement. We are very grateful to the referee for carefully reading the original manuscript and pointing out some typos. The first author would also like to thank Professor Feng Luo, David Gu and Ren Guo for helpful discussions. The research is supported by National Natural Science Foundation of China under grant (No.11501027), and Fundamental Research Funds for the Central Universities (Nos. 2015JBM103, 2014RC028 and 2016JBM071 and 2016JBZ012)..

\bibliographystyle{plain}

Huabin Ge: hbge@bjtu.edu.cn

Department of Mathematics, Beijing Jiaotong University, Beijing 100044, P.R. China\\

Wenshuai Jiang: jiangwenshuai@pku.edu.cn

School of Mathematical Sciences, Zhejiang University, Zheda Road 38, Hangzhou, Zhejiang 310027, P.R. China

\end{document}